\documentstyle{amltd}
\begin{document}
\annalsline{158}{2003}
\received{March 21, 2002}
\startingpage{345}
\def\bye{\end{document}}
 \font\tenrm=cmr10
\def\ritem#1{\item[{\rm #1}]}
\catcode`\@=11
\font\twelvemsb=msbm10 scaled 1100
\font\tenmsb=msbm10
\font\ninemsb=msbm10 scaled 800
\newfam\msbfam
\textfont\msbfam=\twelvemsb  \scriptfont\msbfam=\ninemsb
  \scriptscriptfont\msbfam=\ninemsb
\def\msb@{\hexnumber@\msbfam}
\def\Bbb{\relax\ifmmode\let\next\Bbb@\else
 \def\next{\errmessage{Use \string\Bbb\space only in math
mode}}\fi\next}
\def\Bbb@#1{{\Bbb@@{#1}}}
\def\Bbb@@#1{\fam\msbfam#1}
\catcode`\@=12

 \catcode`\@=11
\font\twelveeuf=eufm10 scaled 1100
\font\teneuf=eufm10
\font\nineeuf=eufm7 scaled 1100
\newfam\euffam
\textfont\euffam=\twelveeuf  \scriptfont\euffam=\teneuf
  \scriptscriptfont\euffam=\nineeuf
\def\euf@{\hexnumber@\euffam}
\def\frak{\relax\ifmmode\let\next\frak@\else
 \def\next{\errmessage{Use \string\frak\space only in math
mode}}\fi\next}
\def\frak@#1{{\frak@@{#1}}}
\def\frak@@#1{\fam\euffam#1}
\catcode`\@=12
\def\mathbb#1{{\Bbb{#1}}}
\newcommand{\p}{\partial}
\newcommand{\e}{e_\alpha}
\newcommand{\sig}{\Sigma}
\newcommand{\x}{\cdot}
\newcommand{\n}{\nabla}
\newcommand{\z}{\bar{z}}

\newcommand{\la}{\langle}
\newcommand{\ra}{\rangle}

\newcommand{\der}[1]{\frac{\partial}{\partial #1}}
\newcommand{\fder}[2]{\frac{\partial #1}{\partial #2}}

\newcommand{\il}[1]{\int\limits_{#1}}

\newcommand{\tf}{\tilde{f}}
\newcommand{\tF}{\tilde{F}}
\newcommand{\tu}{\tilde{u}}
\newcommand{\tg}{\tilde{g}}
\newcommand{\tG}{\tilde{\Gamma}}
\newcommand{\vp}{\varphi}

\newcommand{\mE}{{\cal E}}
\newcommand{\mP}{{\cal P}}
\newcommand{\prn}{P_{u_i}^N}
\newcommand{\prm}{P_{u_i}^{M^\perp}}
\newcommand{\an}{\|_{1,2\alpha}}

\title{Fundamental groups of manifolds\\ with positive isotropic curvature}
\shorttitle{Fundamental groups of PIC manifolds}  

 \acknowledgements{The author was partially supported by NSF Grant DMS-9971927.}
 \author{Ailana M. Fraser}
\institutions{Brown University,
        Providence, RI\\
{\eightpoint {\it Current address\/}}: The University of British Columbia, Vancouver,
BC, Canada\\
{\eightpoint {\it E-mail address\/}:
afraser@math.ubc.ca}}

\centerline{\bf Abstract}
\vglue12pt

A central theme in Riemannian geometry is understanding the relationships 
between the curvature and the topology of a Riemannian manifold. Positive
isotropic curvature (PIC) is a natural and much studied curvature condition
which includes manifolds with pointwise quarter-pinched sectional curvatures 
and manifolds with positive curvature operator. By the results of Micallef and 
Moore there is only one topological type of compact simply connected manifold
with PIC; namely any such manifold must be homeomorphic to the sphere.
On the other hand, there is a large class of nonsimply connected manifolds
with PIC. An important open problem has been to understand the
fundamental groups of manifolds with PIC. In this paper we prove a new
result in this direction. We show that the fundamental group of
a compact manifold $M^n$ with PIC, $n \geq 5$, does not contain a subgroup 
isomorphic to $\mathbb{Z} \oplus \mathbb{Z}$. 
The techniques used involve minimal surfaces.

\section{Introduction}

It is a fundamental problem in geometry to determine the relationships
between the curvature and topology of manifolds. In this paper we study 
the fundamental groups of compact manifolds with positive isotropic 
curvature (PIC). PIC is a natural and much studied curvature condition, 
which first derived its importance from the following beautiful theorem 
of Micallef and Moore~\cite{MM}:

\proclaimtitle{Micallef-Moore}
\proclaim{Theorem}  Let $M$ be a compact $n$\/{\rm -}\/dimensional 
Riemannian manifold with {\rm PIC,} $n \geq 4$. 
Then $\pi_k(M)=0$ for $k=2, \ldots , [\frac{n}{2}]$ {\rm (}\/where $[\cdot]$ denotes
the integer part of the number\/{\rm ).} In particular{\rm ,} if $M$ is simply
connected{\rm ,} then $M$ is homeomorphic to a sphere.
\endproclaim 

If $M$ is an $n$-dimensional Riemannian manifold, one can consider the 
complexification $T_pM \otimes \mathbb{C}$ of the tangent space $T_pM$ 
at the point $p$. The inner product on the tangent space $T_pM$
can be extended to the complexified tangent space 
$T_pM \otimes \mathbb{C}$ as a complex bilinear form 
$( \cdot , \cdot )$ or as a Hermitian inner product 
$\la \cdot , \cdot \ra$. The relationship between these extensions
is given by 
$\la v,w \ra = (v , \bar{w})$ for $v$, $w \in T_pM \otimes \mathbb{C}$.
The curvature tensor extends to complex vectors by linearity,
and the {\it complex sectional curvature} of a two-dimensional subspace 
$\pi$ of $T_pM \otimes \mathbb{C}$ is defined by
$K(\pi)=\la R(v,w)\bar{w}, v \ra$, where $\{v , w \}$ is any
unitary basis of $\pi$. A subspace $\pi$ is said to be {\it isotropic}
if every vector $v \in \pi$ has square zero; that is, 
$(v,v)=0$. 
\vglue4pt {\it Definition} 1.2.
A Riemannian manifold $M$ has {\it positive isotropic curvature}
(PIC) if $K(\pi) > 0$ for every isotropic two-plane 
$\pi$.
\vglue4pt\advance\theoremcount by 1

This curvature condition is nonvacuous only for $n \geq 4$,
since in dimensions less than four there are no two-dimensional
isotropic subspaces.
The classical conditions of pointwise quarter-pinched sectional
curvatures and positive curvature operator are easily seen to
imply PIC; and so in particular, the Micallef-Moore theorem 
gives a generalization of the classical sphere theorem.
Just as positive sectional curvature
is well suited to studying the stability of geodesics, positive
isotropic curvature is in a similar way ideally suited to studying the 
stability of minimal surfaces. The proof of Micallef-Moore involves 
an amplified version of the celebrated existence
theory for minimal two-spheres of Sacks-Uhlenbeck [SU]
together with estimates giving lower bounds on the Morse index
for the area function of any minimal two-sphere in a manifold
with PIC.

By the result of Micallef and Moore, there is only one topological type
of compact simply connected manifold with PIC. On the
other hand there is a large class of nonsimply connected manifolds with
PIC, the simplest being ${S}^1\times {S}^{n-1}$.
This example, with infinite fundamental group $\mathbb{Z}$,
shows in particular that PIC does not imply positive Ricci curvature.
PIC does however imply positive scalar curvature [MW].
It was shown by Micallef and M. Wang [MW] 
that the connected sum of manifolds of PIC also carries such a metric. 
Thus, the fundamental group of a manifold with PIC can be very large.
It is conjectured that
{\it the fundamental group of a manifold with {\rm PIC} is a virtually
free group {\rm (}\/that is{\rm ,} contains a free subgroup of finite index}\/).
Evidence for this is suggested by classification results
of Hamilton~\cite{H}, Micallef-Wang~\cite{MW}
and Noronha~\cite{N1}, ~\cite{N2} 
for four-manifolds with PIC.
For example, an important property of PIC is that it is preserved under the
Ricci flow, and
for $n=4$ Hamilton used the Ricci flow to give a classification
of compact four-manifolds with PIC with no essential incompressible
space forms.
In fact, he has shown that such manifolds 
are diffeomorphic to those which
can be obtained from ${S}^4$, ${S}^1\times {S}^3$ and
${R}P^4$ by the connected sum operation.
Also note that for 
four-dimensional conformally flat manifolds, PIC is equivalent
to positive scalar curvature, 
and it follows from results of
Schoen and Yau~\cite{SY2} that the conjecture holds in this case.
There are a few results in higher dimensions. 
Noronha and Mercuri [MN] for example, show that
the fundamental group of a compact hypersurface $M$ in $\mathbb{R}^k$
with nonnegative isotropic curvature is a free group on 
$b_1(M)$ elements, where $b_1(M)$ denotes the first Betti number
of $M$. 
Another result in \cite{MW} 
is that the second Betti number
of any closed even-dimensional manifold of PIC vanishes, $b_2(M)=0$  
(in particular, a K\"{a}hler manifold can never have PIC).
Little else has been known in dimensions greater than four.
The main result of this paper is:
\proclaim{Theorem} \label{theorem:pic}
Let $M$ be a compact $n$\/{\rm -}\/dimensional Riemannian manifold with positive 
isotropic curvature{\rm ,} $n \geq 5$.
Then the fundamental group of $M$ does not contain a subgroup isomorphic
to $\mathbb{Z} \oplus \mathbb{Z}$.
\endproclaim
  
This proves a weaker version of the above conjecture, that
any {\it free} abelian subgroup of the fundamental group is 
cyclic.
This may be thought of as a PIC version of the Chern conjecture for manifolds
with positive sectional curvature (recently proven false by 
Shankar~\cite{Sh}); 
we note that the assumption that the
abelian subgroup be {\it free} cannot be removed under the PIC condition 
(quotients of ${S}^1\times {S}^{n-1}$).

The idea of the proof of Theorem \ref{theorem:pic} is as follows. 
If the fundamental group did contain a subgroup isomorphic to
$\mathbb{Z} \oplus \mathbb{Z}$, then by results of Schoen and Yau~\cite{SY1}
it could be represented by an area minimizing torus $T^2$.
There are manifolds with PIC that do have stable tori:

\numbereddemo{{E}xample}
$S^1 \times L$, where $L$ is a lens space, admits a metric of PIC.
We can define a map $h:T^2 \rightarrow S^1 \times L$ such that 
$h_*: \pi_1(T^2) \rightarrow \pi_1(S^1\times L)=\mathbb{Z}\oplus\mathbb{Z}_p$
is onto by viewing $T^2$ as the square with sides identified, mapping the 
boundary to
generators of $\pi_1(S^1 \times L)$ and extending to the interior via
contractibility. Any simple closed curve on $T^2$ has the form 
$\pm ra \pm sb$ where $a$ and $b$ represent generators of
$\pi_1(T^2)$ and $r$ and $s$ are relatively prime integers,$(r,s)=1$. Assume
$a$ and $b$ are mapped to generators $\alpha$ and $\beta$ of 
$\pi_1(S^1 \times L)$. Then the image of the simple closed curve is
$\pm r\alpha \pm s\beta$, which is nontrivial in 
$\pi_1(S^1 \times L)$ since $(r,s)=1$. Thus, $u$ maps any simple closed
curve in $T^2$ to a noncontractible curve in $S^1 \times L$, and
it then follows by the proof of \cite{SY1} that there is a branched 
minimal immersion $u:T^2 \rightarrow S^1 \times L$ such that $u_*=h_*$ on
$\pi_1(T^2)$ and the induced area of $u$ is least among all maps with
the same action on $\pi_1(T^2)$.  
\enddemo

However, we argue that for sufficiently large $k$, any minimal torus
representing a subgroup
$k \mathbb{Z} \oplus k \mathbb{Z}$ of the fundamental group,
must be unstable.
To prove instability we use a complex formula for the second
variation of area developed by Siu and Yau~\cite{SiY} in their proof of 
the Frankel
conjecture, and by Micallef in~\cite{Mi}, 
and Micallef and Moore~\cite{MM}.
PIC is naturally suited to studying the stability of minimal surfaces 
because one of the terms in the complexified second variation
formula involves isotropic curvature; and, in order to prove instability
or estimate the index of a minimal surface one must find 
holomorphic isotropic deformations of the surface.
The argument for the index estimate for minimal 
two-spheres in manifolds with PIC of Micallef and Moore uses
Grothendieck's splitting theorem for holomorphic
vector bundles over the Riemann sphere.
In the case of a torus, there is not such a nice splitting theorem,
and in fact holomorphic deformations will not exist in general.
However, for large $k$ we are able to find {\it almost holomorphic}
isotropic deformations that violate the stability inequality.
We expect to be able to extend our techniques from the
torus to the case of higher genus surfaces.
In fact the techniques seem promising for proving
the full conjecture.

Although our techniques are different, the idea of finding almost
holomorphic sections is reminiscent of Donaldson's~\cite{D}
almost holomorphic sections used in his procedure for producing 
symplectic submanifolds of any even
codimension within a symplectic manifold.
Also, in finding our sections, we use the fact that for 
sufficiently large $k$, there is a distance decreasing,
degree one mapping from the torus to the two-sphere.
The existence of such a map is somewhat analogous to 
Gromov and Lawson's~\cite{GL} notion of `enlargeability',
which was used by them in studying manifolds with positive scalar
curvature.
\vglue5pt
{\it Acknowledgement}.\ I would like to thank Richard Schoen for many
helpful and inspiring conversations throughout the progress of this work.

\vglue-8pt
\section{Proof of the theorem}
\label{section:one}
\vglue-4pt

In this section we prove the main result, Theorem \ref{theorem:pic}.
Let $(M, g)$ be a compact $n$-dimensional Riemannian manifold with PIC, 
$n \geq 5$, and assume that $\pi_1(M)$ contains a subgroup isomorphic to 
$\mathbb{Z} \oplus \mathbb{Z}$.
Let $\kappa >0$ denote a lower positive bound on the isotropic
curvature of $M$.
Choose $\varepsilon$ such that $0 < 2(C\varepsilon)^2 < \kappa$ where
$C$ is to be chosen later.
Given an integer $k \geq 1$, consider the subgroup 
$G=k\mathbb{Z} \oplus k\mathbb{Z}$ of $\pi_1(M)$.
By the results of Schoen-Yau~\cite{SY1} (and Sacks-Uhlenbeck~\cite{SU2})
there is a conformal branched minimal immersion of the torus 
$u: T^2 \rightarrow M$ whose induced map $u_{*}$
on the fundamental group has image $G$, with least area among branched 
immersions with the same action as $u_{*}$. However, we claim
that for $k$ sufficiently large, such a minimal torus must be
{\it unstable}. A key point in proving this is the following:
For $k$ sufficiently large, there is a distance decreasing, degree one
map $f: (T^2, u^*g) \rightarrow (S^2, h)$, where $h$ is the standard
metric on $S^2$, with $|df| < \varepsilon$.
The existence of this map follows from Lemma \ref{lemma:map}
from the Appendix. From now on, we fix $k$ so that the lemma holds.

Let $\Sigma=u(T^2)$ and let $N\Sigma$ be the normal bundle to the
surface $\Sigma$ in $M$. Consider the pull-back of the 
normal bundle $u^*(N\Sigma)$ with the pull back of the metric and
normal connection $\nabla^{\perp}$.
Let $E=u^*(N\Sigma) \otimes \mathbb{C}$ be the complexified bundle.
The metric on $u^*(N\Sigma)$ extends as a complex bilinear form
$( \cdot , \cdot )$ or as a Hermitian metric $\la \cdot , \cdot \ra$
on $E$, and the connection $\n^{\perp}$ and curvature tensor extend
complex linearly to sections of $E$.
There is a unique holomorphic structure on $E$ such that the
$\bar{\p}$ operator
$$
      \bar{\p}: {\cal A}^{p,q}(E) \rightarrow {\cal A}^{p,q+1}(E),
$$
where ${\cal A}^{p,q}(E)$ denotes the space of $(p,q)$-forms on
$T^2$ with values in $E$, is given by
$$
      \bar{\p}\omega=(\nabla_{\der{\bar{z}}}^{\perp} \omega)d\bar{z}
$$
where $\der{\bar{z}}=\frac{1}{2}(\der{x}+i\der{y})$, for local coordinates
$x, y$ on $T^2$. 
Suppose $\Sigma$ is {\it stable}. Then the complexified
stability inequality (see \cite{SiY}, \cite{Mi}, \cite{MM}, \cite{F}) holds:
$$
    \il{T^2} \la R(s,\fder{u}{{z}})\fder{u}{\bar{z}},s \ra \; dxdy
    \leq \il{T^2} [|\nabla_{\der{\bar{z}}}^{\perp} s|^2 
    - |\nabla_{\der{z}}^{\top} s|^2] \; dxdy
$$
for all $s \in \Gamma(E)$. Assume now that $s$ is isotropic. Since $u$
is conformal, $\fder{u}{z}$ is isotropic and $\{s,\fder{u}{z}\}$ span
an isotropic two-plane. Using the lower bound on the isotropic
curvature and throwing away the second term on the right, we get
\begin{equation} \label{equation:stability}
    \kappa \il{T^2} |s|^2 \; da 
    \leq \il{T^2} |\bar{\p} s|^2 \; da
\end{equation}
where $da$ denotes the area element for the induced metric $u^*g$ on $T^2$.
We now argue that we can find an ``almost holomorphic'' isotropic
section of $E$ that violates this stability inequality.
That is, we will find $s \in \Gamma(E)$ such that
$$
     \il{T^2} |\bar{\p} s|^2 \; da 
    << \il{T^2} |s|^2 \; da.
$$

Recall from \cite[p.\  209]{MM}, that the bilinear pairing
$( \cdot , \cdot )$ establishes a holomorphic isomorphism between
$E$ and its dual $E^*$, and hence $c_1(E)=0$.
Let $L$ be a positive holomorphic line bundle
over $S^2$ with metric and connection, with $c_1(L) \geq 2$. 
Pull back $L$ to a line bundle
$\xi= f^*L$ over $T^2$ with pull-back metric and connection. The
pull-back metric and connection on $\xi$ determine a holomorphic
structure on $\xi$, and
$$
    c_1(\xi)=(\deg f)c_1(L)=c_1(L) \geq 2.
$$
Now consider the tensor product holomorphic vector bundle 
$\xi \otimes E$ of rank $n-2$ over $T^2$.
Let ${\cal H}$ denote the complex vector space of holomorphic
sections of $\xi \otimes E$.
By the Riemann-Roch theorem for holomorphic vector bundles over a 
Riemann surface (\cite[p.\ 64]{Gu}) 
we have
\begin{eqnarray} 
     \dim {\cal H}(\xi \otimes E) 
     & \geq&  c_1( \det(\xi \otimes E) ) \\
     & =&  \hbox{rank}(E) c_1(\xi) + c_1(E) \nonumber\\
     & = & (n-2) c_1(\xi) .\nonumber
\label{equation:dim}
\end{eqnarray}
Therefore, there is at least a $2n-4$ dimensional space of holomorphic
sections of the tensor product bundle $\xi \otimes E$.
We now prove existence of a holomorphic section of
$\xi \otimes E$ which is {\it isotropic}.
Define a complex bilinear pairing
$$
    {\cal H}(\xi \otimes E) \times {\cal H}(\xi \otimes E)
    \rightarrow {\cal H}(\xi \otimes \xi)
$$
given by
$$
    (t \otimes s_1 , t \otimes s_2)
    = t \otimes t \; (s_1 , s_2)
$$
on any local basis of ${\cal H}(\xi \otimes E)$ of the form
$\{t \otimes s_i \}_{i=1}^{n-2}$, where $t$, $s_i$ ($i=1, \ldots,\break n-2$)
are local holomorphic sections giving a local basis for holomorphic 
sections of $\xi$ and $E$
respectively.
Given any $x \in T^2$ we obtain a homogeneous polynomial on 
$\mathbb{C}^m \cong {\cal H}(\xi \otimes E)$ where 
$m=\dim {\cal H}(\xi \otimes E)$, given by
$P_x(\sigma)=(\sigma, \sigma)(x)$.
The zero set
$$
   V(P_x)=\{ \sigma \in \mathbb{P}^{m-1}: P_x(\sigma)=0\}
$$
is a hypersurface in $\mathbb{P}^{m-1}$.
Now given $m-1$ distinct points, we obtain $m-1$ hypersurfaces
in $\mathbb{P}^{m-1}$, and observe that $m-1$ hypersurfaces
in $\mathbb{P}^{m-1}$ intersect in a nonempty set of points.
The intersection is a set of holomorphic sections of 
$\xi \otimes E$ which are isotropic at $m-1$ distinct points.
Let $\sigma \in {\cal H}(\xi \otimes E)$ be such a section. 
Then $(\sigma, \sigma)$ is a holomorphic section of 
$\xi \otimes \xi$ with at least $m-1$ zeros. But the number
of zeros of a holomorphic section of $\xi \otimes \xi$ is
$2c_1(\xi)$. From (2.2)
$$
      m-1 \geq (n-2)c_1(\xi)-1 > 2c_1(\xi)
$$
since $n \geq 5$ and we chose $c_1(\xi) \geq 2$.
It follows that $(\sigma, \sigma) \equiv 0$ and so $\sigma$ is
isotropic.

Finally, we show that any holomorphic isotropic section of 
$\xi \otimes E$ produces an {\it almost} holomorphic isotropic section of $E$.
Let $\tilde{s} \in \Gamma(\xi \otimes E)$ be holomorphic and
isotropic.
First, we choose $t_1^*$, $t_2^* \in \Gamma(L^*)$ that satisfy
$$
      |t_1^*| + |t_2^*| \geq 1
$$
and
$$
      |t_1^*| + |t_2^*| \leq 2.
$$
This is possible since one can for example let
$t_1^*$ be a trivialization of the dual bundle $L^*$ over
$S^2-U_-$, where $U_-$ is some neighborhood of the south pole in
the southern hemisphere $S^2_-$ of $S^2$, with pointwise norm
$|t_1^*|=1$. Now extend $t_1^*$ to a global
section of $L^*$ and multiply by a cut-off function 
$\varphi$, $0 \leq \varphi \leq 1$ with $\varphi \equiv 1$ on $S^2_+$ 
(the northern hemisphere of $S^2$) and $\varphi \equiv 0$ on $U_-$.
Similarly, we obtain $t_2^* \in \Gamma(L^*)$ with 
$|t_2^*| \leq 1$ on $S^2$ and $|t_2^*| = 1$ on $S^2_-$.

Now pull back $t_1^*$, $t_2^*$ to 
$\alpha_1=f^*t_1^*$, $\alpha_2=f^*t_2^* \in \Gamma(\xi^*)$.
We can think of $\alpha_i$ ($i=1,2$) as a linear mapping
$$
     \alpha_i : \Gamma(\xi \otimes E)
     \rightarrow \Gamma(E)
$$
where $\alpha_i$ is defined on any local basis section 
$t \otimes s_j$ by 
$$
     \alpha_i(t \otimes s_j) = \alpha_i(t) s_j.
$$
Note that the image of an isotropic section of $\xi \otimes E$
under $\alpha_i$ is an isotropic section of $E$.
Define 
$s_1=\alpha_1(\tilde{s})$, $s_2=\alpha_2(\tilde{s}) \in \Gamma(E)$
and observe that 
\begin{eqnarray*}
     |s_1|^2 + |s_2|^2 
     & = & |\alpha_1(\tilde{s})|^2 + |\alpha_2(\tilde{s})|^2  \\
     & = & |\alpha_1|^2|\tilde{s}|^2 + |\alpha_2|^2|\tilde{s}|^2 \\
     & = & ( |\alpha_1|^2 + |\alpha_2|^2) |\tilde{s}|^2 \\
     & \geq & |\tilde{s}|^2.
\end{eqnarray*}
Integrating this pointwise inequality 
we obtain
$$
    \il{T^2} |s_1|^2 \;da + \il{T^2} |s_2|^2 \;da
    \geq \il{T^2} |\tilde{s}|^2 \;da.
$$
Therefore, for either $i=1$ or $i=2$, say $i=1$, we must have
\begin{equation} \label{equation:section}
    \il{T^2} |s_1|^2 \;da 
    \geq \frac{1}{2} \il{T^2} |\tilde{s}|^2 \;da.
\end{equation}
Now,
\begin{eqnarray*}
       |\bar{\p} s_1| & = & | \bar{\p} (\alpha_1(\tilde{s})) | \\
    & = & | (\bar{\p} \alpha_1) \tilde{s} + \alpha_1(\bar{\p}\tilde{s})| \\  
    & = & | (\bar{\p}(f^*t_1^*) \tilde{s}| \\
    & \leq & C |\bar{\p}f||\tilde{s}| \\
    & \leq & C \varepsilon |\tilde{s}|
\end{eqnarray*}
where $C$ depends only on $L$.
Combining this with (\ref{equation:section}) we have
$$
    \il{T^2} |\bar{\p} s_1|^2 \; da
    \leq 2(C\varepsilon)^2 \il{T^2} |s_1|^2 \; da
    < \kappa \il{T^2} |s_1|^2 \; da.
$$
Therefore $s_1$ is an {\it almost holomorphic} isotropic section
of $E$, that violates the stability inequality (\ref{equation:stability}).
This completes the proof of Theorem \ref{theorem:pic}.
\\

By a similar argument we also obtain the following.
\proclaim{Theorem}
There exists a sufficiently high finite cover of any stable incompressible
minimal torus in a manifold with positive isotropic curvature which is 
unstable.
\endproclaim

\section{Appendix: Distance decreasing map}
\label{section:map}

\proclaim{Lemma} \label{lemma:map}
Let $M$ be a compact Riemannian manifold whose fundamental group contains a 
subgroup isomorphic to $\mathbb{Z} \oplus \mathbb{Z}$. 
Given $\varepsilon >0${\rm ,} there exists $k$ sufficiently large such that for
any torus $\Sigma$ in $M$
representing the subgroup $G=k\mathbb{Z} \oplus k\mathbb{Z}$ 
of $\pi_1(M)${\rm ,} there is a distance
decreasing{\rm ,} degree one mapping $f: \Sigma \rightarrow S^2$ with 
$|df| < \varepsilon$. 
\endproclaim

\demo{Proof}
Let $\tilde{M}$ denote the universal cover of $M$.
First observe that for $k$ large, 
any torus $\Sigma$ in $M$ representing the group 
$G= k\mathbb{Z} \oplus k\mathbb{Z}$
has large systole. Recall that the systole (see \cite{G}) of
$\Sigma$ is defined to be the number
$$
      {\cal L}=\inf \{ l(\gamma): 
      \gamma\; \hbox{a noncontractible closed curve in}\; \Sigma \}.
$$
To see this, define the systole associated to the group element 
$\gamma \in \pi_1(M)$,
$$  
        \delta(\gamma)=\inf_{x \in \tilde{M}} d(x, \gamma x);
$$
that is, $\delta(\gamma)$ is the length of the shortest  closed curve in
$M$ freely homotopic to $\gamma$. Note that since $M$ is compact, the set
$S=\{\gamma \in \pi_1(M): \delta(\gamma) < C \}$ is finite for any fixed 
constant
$C$. Thus, for $k$ large, $G \cap S=\phi$. That is, given $L$, there exists
$k$ such that the systole of $\Sigma$ is greater than $L$. 

Let $\tilde{\Sigma}$ be the universal cover of $\Sigma$, and
$p:\tilde{\Sigma} \rightarrow \Sigma$ the covering map.
Since $\tilde{\Sigma}$ is complete and noncompact with compact quotient
$\Sigma$, there is a geodesic line 
$r: \mathbb{R} \rightarrow \tilde{\sig}$.
Let ${\cal C}_1$ be the component of $\tilde{\sig}-r(\mathbb{R})$
which is on the side of $\tilde{\sig}$ in the direction of the unit 
normal $\nu$ to $r$ such that $\{r'(0), \nu(0) \}$ is positively oriented.
Choose $T$ very large, $T>>L$.
Define $D_1: \tilde{\sig} \rightarrow \mathbb{R}$ by,
$$
     D_1(x)=d(x,r(T))-T
$$
and let $D_2: \tilde{\sig} \rightarrow \mathbb{R}$ be the signed
distance function to $r$,
$$
   D_2(x)= \left\{ \begin{array}{ll}
    d(x,r)& \hbox{for}\;x \in {\cal C}_1   \\
    -d(x,r)& \hbox{for}\;x \in \tilde{\sig} - {\cal C}_1 .
    \end{array} \right. 
$$
Both $D_1$ and $D_2$ are Lipschitz continuous with derivative
bounded by 1. 
Consider the region
$$
   \tilde{{\cal R}}=\left\{x \in \tilde{\Sigma} : 
   |D_1(x)| \leq \frac{L}{4},
   |D_2(x)| \leq \frac{L}{4} \right\}.
$$
Define the map $f: \tilde{{\cal R}} \rightarrow \mathbb{R}^2$ by
$$
     f(x)= (D_1(x),D_2(x)).
$$
Then the boundary of $\tilde{{\cal R}}$ is mapped to the boundary of the 
rectangle
$[-\frac{L}{4},\frac{L}{4}]\times[-\frac{L}{4},\frac{L}{4}]$ 
in $\mathbb{R}^2$.
Also, $r(0)$ is the only point in $\tilde{{\cal R}}$ which is mapped under 
$f$ to the origin in $\mathbb{R}^2$. 
In fact, $f$ is a local diffeomorphism in a neighborhood of $r(0)$, 
and hence the degree of $f$ is equal to one on the component of 
$\mathbb{R}^2-f(\p \tilde{{\cal R}})$ containing the origin.
In the definition of $f$, now replace $D_1$, $D_2$ by smooth approximations
of these functions, which we will denote by $\bar{D}_1$ and $\bar{D}_2$,
with $\|D_i-\bar{D}_i\|_{C^1} <\delta$ for $i=1,2$,
where $\delta$ may be chosen arbitrarily small. For this smoothly 
redefined $f$ also, by the above argument, 
$f(\tilde{{\cal R}})$ covers a disk of radius at least 
$\frac{L}{5}$ about the origin in $\mathbb{R}^2$.

Let $\lambda : D(0,\frac{L}{5}) \rightarrow D(0,\pi)$ be the
contraction $\lambda(x)=\frac{5\pi}{L}x$ (where $D(0,s)$ denotes
the disk of radius $s$ centered at the origin in $\mathbb{R}^2$). Let 
$e: D(0,\pi) \rightarrow S^2$ be the exponential map at the north
pole $n$ of the sphere $S^2$ with the standard metric, 
$e(x)=\exp_n(x)$. Extend $g=e \circ \lambda$ to $f(\tilde{{\cal R}})$
by defining $g \equiv s$ (the south pole) on 
$f(\tilde{{\cal R}}) - D(0,\frac{L}{5})$. Then 
$\tilde{F}=g \circ f : \tilde{{\cal R}} \rightarrow S^2$ 
is smooth with derivative
$$
     |d\tilde{F}| \leq |df||d\lambda||de| \leq \frac{C}{L}
$$
where $C$ is a constant (independent of $k$).

Observe that 
diam $\tilde{{\cal R}} \leq L$. Given $x, y \in \tilde{{\cal R}}$,
we have
$d(x,r)=d(x,r(t_1)) \break\leq \frac{L}{4}$ and
$d(y,r)=d(y,r(t_2)) \leq \frac{L}{4}$ 
for some $t_1, t_2$ with $r(t_1), r(t_2) \in \tilde{{\cal R}}$, 
and
$$
    d(r(t_1),r(t_2)) \leq D_1(r(t_1))
                                 +  D_1(r(t_2))
    \leq \frac{L}{2}.
$$
Now,
$$
    d(x,y) \leq d(x, r(t_1)) + d(r(t_1),r(t_2))
                 + d(r(t_2),y) \leq L.
$$
Hence no two points of $\tilde{{\cal R}}$ are identified in the quotient
$\Sigma$. If two points were identified, then any minimal curve in
$\tilde{{\cal R}}$ joining the two points would project to a nontrivial 
closed curve of length less than or equal to $L$, which is less 
than the systole, a contradiction.

Therefore, $\tilde{{\cal R}}$ projects one to one into $\Sigma$, and
we may define $F: {\cal R} \rightarrow S^2$, 
where ${\cal R}=p(\tilde{{\cal R}})$, by
$F(x)=\tilde{F}(p^{-1}(x))$. 
Since $F$ maps $\p {\cal R}$ to $s$, we can extend  
$F$ from ${\cal R}$ to a map $F: \sig \rightarrow S^2$ with the
desired properties, by defining $F \equiv s$ on $\sig - {\cal R}$.
\enddemo

\AuthorRefNames [MSY]

\end{document}